\documentclass[oneeqnum]{article}  
\usepackage{amsmath,amscd,amssymb,epsf}
\usepackage{amsfonts,graphicx}
\usepackage[usenames,dvipsnames]{color} 	
\usepackage{epsfig}
\usepackage{graphics}
\usepackage{multirow}
\usepackage{rotating}  
\usepackage{pdflscape} 
\usepackage{float}
\usepackage{soul} 

\usepackage{verbatim}

\usepackage{algorithm}         
\usepackage{algpseudocode}
\algblock{ParFor}{EndParFor}
\algnewcommand\algorithmicparfor{\textbf{parfor all}}
\algnewcommand\algorithmicpardo{\textbf{do}}
\algnewcommand\algorithmicendparfor{\textbf{end\ parfor}}
\algrenewtext{ParFor}[1]{\algorithmicparfor\ #1\ \algorithmicpardo}
\algrenewtext{EndParFor}{\algorithmicendparfor}

\usepackage{subfloat}
\usepackage{booktabs}
\usepackage{subcaption}

\usepackage{acronym}
\usepackage{txfonts}

\usepackage{setspace}
\usepackage{array}                      

\newcommand{\Texp}{\mathbf{t}^{\mbox{{\tiny \rm exp}}}}
\newcommand{\Texpdif}{\mathbf{t}^{\mbox{{\tiny \rm exp}}}_{\mbox{{\tiny \rm dif}}}}

\newcommand{\dbar}{\bar{\partial}}

\newcommand{\T}{{\mathbf{t}}}
\newcommand{\Om}{\Omega}
\newcommand{\DOm}{\partial\Omega}
\newcommand{\R}{{\mathbb R}}

\newcommand{\TexpR}{\mathbf{t}^{\mbox{{\tiny \rm exp}}}_{\mbox{{\tiny R}}}}

\newcommand{\muEXPR}{\mu^{\mbox{{\tiny exp}}}_{\mbox{{\tiny R}}}}

\begin{document}

\title{A Real-time D-bar Algorithm for 2-D Electrical Impedance Tomography Data}

\author{Melody Dodd \thanks{Department of Mathematics, Colorado State University, USA}, Jennifer L. Mueller\thanks{Department of Mathematics and School of Biomedical Engineering, Colorado State University, USA}}

\date{}

\maketitle

\begin{abstract} The aim of this paper is to show the feasibility of the D-bar method for real-time 2-D EIT reconstructions.  A fast implementation of the D-bar method for reconstructing conductivity changes on a 2-D chest-shaped domain is described. Cross-sectional difference images from the chest of a healthy human subject are presented, demonstrating what can be achieved in real time.  The images constitute the first D-bar images from EIT data on a human subject collected on a pairwise current injection system.
\end{abstract}

\section{Introduction}
The generalized Laplace equation serves as a model for the electric field $u(x)$ propagating in living tissue from a low frequency, low amplitude applied current density on the surface of a body.  Denoting the conductivity by $\sigma(x)$ and a bounded 2-D domain by $\Om$, the governing equation is
\begin{equation}\label{genLap}
\nabla\cdot( \sigma(x)\nabla u(x))=0, \quad x\in\Om.
\end{equation}
The ideal data is the Dirichlet-to-Neumann (DN) map, or voltage-to-current density map $\Lambda_{\sigma}$, defined by
\begin{equation} \label{DNmap}
\Lambda_{\sigma}:u\vert_{\DOm}\rightarrow \sigma(x)\frac{\partial u}{\partial\nu}\big|_{\DOm},
\end{equation}
where $u$ is a solution to \eqref{genLap} and $\nu$ is the unit normal vector to the surface.  Equation \eqref{genLap} serves as the governing equation for the electric field in electrical impedance tomography (EIT), and has a rich mathematical history dating back to the problems posed by Calder\'on \cite{Calderon}:  (1)  when does the inverse problem of determining $\sigma$ from knowledge $\Lambda_{\sigma}$  have a unique solution, and (2) how can it be determined?  Historical reviews of the answers to these questions can be found in \cite{borcea_review, JenSamuBook}, and the reader will find that most of the uniqueness proofs have utilized complex geometrical optics (CGO) solutions.  Some have also been formulated as constructive proofs \cite{nachman, brownuhlmann, astala}, and most of these include PDEs known as D-bar, or $\dbar$, equations for the CGO solutions.  These works have given rise to a new class of direct (noniterative) reconstruction algorithms known as D-bar methods.

 D-bar equations are of the form $\dbar w = f,$ where $f$ may depend on $w$, and the $\dbar$ operator is defined by $\dbar_z = \frac{1}{2} \left(\partial_x +i \partial_y\right)$, where $z=x+iy$.  The common threads in these methods are (1) a direct relationship between the CGO solutions and the unknown conductivity, (2) a nonlinear Fourier transform, also known as the scattering transform, providing a link between the data and the CGO solution, and (3) a D-bar equation to be solved for the CGO solutions with respect to a complex-frequency variable. 
These steps can be expressed in general terms by
$$\mbox{DN map} \longrightarrow \mbox{scattering transform} \longrightarrow \mbox{CGO solution} \longrightarrow \mbox{conductivity}$$

To briefly review the history of D-bar methods for the inverse conductivity problem in dimension 2, we begin with the constructive global uniqueness proof  for twice differentiable conductivities \cite{nachman}.  The fast method implemented in this paper is based on that work and subsequent theory and implementations including \cite{siltanen2000,properties,siam,TMIdata,ChestPaper,murphy2009,deangelo,dbar_regul}.   The regularity conditions on the conductivity in \cite{nachman} were relaxed to once-differentiable in  \cite{brownuhlmann}, and a reconstruction method based on this work can be found in \cite{knudsen,knudsenMummy,knudsenTamasan}.  A nearly constructive proof for complex conductivities of the form $\sigma+i\omega \epsilon$ with  $\sigma,\; \epsilon
 \in W^{2,\infty}(\Omega)$, with $\omega$ small was presented in \cite{francini}, and a direct algorithm and implementations can be found in \cite{hamilton2013,hamilton2012,herrera2014}.  A non-constructive proof that applies to complex admittivities with no smallness assumption is found in \cite{Bukhgeim2007}.  Astala and P\"aiv\"arinta provide a {\sc CGO}-based constructive proof for real conductivities $\sigma \in L^{\infty }(\Omega ),$ and numerical results related to this work can be found in \cite{APpaper1,APpaper2}.

Electrical Impedance Tomography holds great promise as a bedside imaging tool for patients in intensive care. In the 2-D geometry, EIT is clinically useful for chest imaging.  Conductivity images have been used for monitoring pulmonary perfusion \cite{BrownBarber, frerichs_perfusion,Smit2004}, determining  regional ventilation in the lungs \cite{frerichs2007, frerichs2002, victorino},   detecting extravascular lung water \cite{kunst2}, and evaluating shifts in lung fluid in congestive heart failure patients \cite{freimark}.   Regional results have been validated with CT images  \cite{frerichs_perfusion,frerichs2007,Costa2009,Smit2004} and  radionuclide scanning \cite{kunst}  in the presence of pathologies such as atelectasis, pleural effusion, and pneumothorax.  However, the solution of the inverse problem in real-time poses a significant challenge. D-bar methods have been generally regarded as computationally intensive, but in this work, we show that through parallelization and careful optimization of the computational routines,  a fast implementation is capable of providing real-time images from the pairwise current injection system at CSU.

In this work, we chose to optimize the D-bar method based on the uniqueness proof \cite{nachman} and subsequent results and implementations \cite{siltanen2000, murphy2009}.  Many features of the fast implementation also apply to numerical solution methods of other D-bar reconstruction algorithms.  

The paper is organized as follows.  Section \ref{sec:background} contains a brief mathematical description of the D-bar method implemented here.  Section \ref{sec:fast} describes the fast implementation,  parallelized in two different ways.  Section \ref{sec:results} contains tables of runtimes and reconstructions on three different meshes from a set of data collected on a human subject.  The final two sections contain conclusions and acknowledgments.

\section{Background}  \label{sec:background}
We begin with an overview of the D-bar method implemented here, both for the reader's convenience and to place the fast implementation in its mathematical context.  For further details, see \cite{nachman, JenSamuBook}. 

The method begins with a transformation of the generalized Laplace equation with conductivity $\sigma \in W^{2,p}(\Om)$, for some $p>1$, to the 
Schr\"odinger equation through the change of variables $q(z)=\Delta\sqrt{ \sigma(z)}/\sqrt{ \sigma(z)}$ and $\tilde{u}(z)= \sqrt{\sigma(z)}u(z)$, where the point $z=(x,y)$ lies in $\Om$, a bounded, simply connected Lipschitz domain in $\R^2$.  This results in
  \begin{equation}\label{Schr}
 -\Delta \tilde{u}+q(z)\tilde{u}=0, \quad z\in\Om.
 \end{equation}
Under the assumption that $ \sigma$ is constant in a neighborhood of the boundary of $\Om$, one can extend \eqref{Schr} to the whole plane, taking $q=0$ outside $\Om$.  Without loss of generality, we will assume $ \sigma\equiv 1$ in a neighborhood of the boundary.  

The existence of CGO solutions to \eqref{Schr} in the plane was established by Faddeev \cite{faddeev} in the context of quantum physics and shown by Nachman \cite{nachman} to always exist for $q$ of the form $q(z)=\Delta\sqrt{ \sigma(z)}/\sqrt{ \sigma(z)}$.   Introducing the complex parameter $k=k_1+ik_2$, and identifying the spatial variable $z = x+iy$ with the corresponding point in the complex plane, the CGO solution $\psi(z,k)$ satisfies \cite{nachman}
 \begin{eqnarray}\label{Schr_psi}
 -\Delta \psi(z,k)+q(z)\psi(z,k)&=&0, \quad z\in\R^2 \\
e^{-ikz}\psi(z,k) -1 &\in& W^{1,p}(\R^2), \quad p>2. \label{psi_asy}
 \end{eqnarray}
The CGO solution closely related to $\psi$ is $\mu(z,k)$, defined by
$\mu(z,k)\equiv e^{-ikz}\psi(z,k)$.
The conductivity can be obtained directly from $\mu$ or $\psi$ through the formula \cite{nachman}
\begin{equation} \label{sigmaeval}
   \sigma(z)=\mu^2(z,0), \quad z\in\Omega,
\end{equation}
and so the key steps in the method are to link the data $\Lambda_{\sigma}$ to the function $\mu(z,k)$.

These links are provided by the scattering transform $\T(k)$ and the D-bar equation for $\mu$.  The scattering transform is defined by 
\begin{equation}\label{scat_def}
\T(k) = \int_{\Om}e^{i\bar{k}\bar{z}}q(z)\psi(z,k) dz,
\end{equation}
and can be regarded as a nonlinear Fourier transform of $q$ in light of the asymptotic behavior of $\psi$.    The D-bar equation satisfied by $\mu$ is
\begin{equation}\label{Dbar_eq}
  \frac{\partial \mu}{\partial\bar{k}}= \frac{\T(k)}{4\pi\overline{k}}e_{-z}(k)\overline{\mu(z,k)},
\end{equation}
where $e_{z}(k)\equiv e^{i(kz+\bar{k}\bar{z})}$.  The scattering transform is related to the DN data through an equation that requires knowledge of $\psi$ on the boundary of $\Om$:
 \begin{equation}  \label{TkDN}
  \T(k) = \int_{\partial\Om}e^{i\overline{k}\bar{z}}(\Lambda_{ \sigma}-\Lambda_{1})\psi(z,k) ds.
\end{equation}
Here $\Lambda_1$ denotes the DN map corresponding to the homogeneous conductivity distribution $\sigma\equiv 1$.  For the fast implementation, we utilize a linearized approximation to the scattering transform, denoted by $\Texp$, which is defined by replacing $\psi\vert_{\DOm}$ by its asymptotic behavior:
\begin{equation}  \label{Texp}
  \Texp(k) \equiv \int_{\partial\Om}e^{i\overline{k}\bar{z}}(\Lambda_{ \sigma}-\Lambda_{1})e^{ikz} ds.
\end{equation}
This approximation was first introduced in \cite{siltanen2000} and was later studied in \cite{properties} where it was shown that the D-bar equation \eqref{Dbar_eq} with $\T(k)$ replaced by $\Texp$ truncated to a disk of radius $R$ in the $k$-plane has a unique solution which is smooth with respect to $z$, and the reconstruction is smooth and stable.  Further, it was shown that no systematic artifacts are introduced when the method with $\Texp$ is applied to piecewise continuous conductivities.  

In summary, the mathematical steps of the algorithm implemented here are (1) compute $\Texp(k)$ on a disk of radius $R$ in the $k$-plane, (2) solve the D-bar equation \eqref{Dbar_eq} for each $z$ in the region of interest on the disk $|k|\leq R$, and (3) compute $\sigma$ from \eqref{sigmaeval}.


\section{Fast implementation}  \label{sec:fast}
A fast implementation in Matlab on a 12 core Mac Pro with two 2.66 GHz 6 core Intel Xeon processors and Matlab's Parallel Computing Toolbox is capable of computing reconstructions at less than the data acquisition rate of 16 frames/s, or  0.0625 s/frame, of the ACE 1 pairwise current injection EIT system at CSU \cite{hardware}.  This demonstrates the feasibility of CGO methods for real-time reconstructions.  

In fact, we consider two options for the parallel computations.  Ideally, in real-time reconstructions, data is collected, demodulated, and fed directly to the reconstruction algorithm, one frame at a time.  In this configuration, only the loop over the $z$-values in the solution of the D-bar equation is trivially parallelizable.  This accounts for over 95\% of the computation time and can be used to obtain real-time reconstructions at a rate of 0.0621 s/frame on 7 cores on a coarse mesh of 562 elements, as reported in Section \ref{sec:results}.  This approach is structured as shown in Algorithm~\ref{DbarAlg_zloop}.

If a time delay of approximately one second is acceptable to the user, the algorithm can be parallelized over the frames.  In this configuration, the data is collected, demodulated, and sent to a buffer from which multiple frames are sent to the reconstruction algorithm in batch.  This approach results in an algorithm that is over 99\% parallelizable, and the frame rate is even faster.  As shown in Section \ref{sec:results}, reconstructions were computed at a rate of 0.0215 s/frame on 64 cores or 0.0578 s/frame on 12 cores on a mesh of 1931 elements. The computational method for this approach is shown in Algorithm~\ref{DbarAlg_frames}. 


\begin{algorithm}[h]
\caption{Fast Parallelized D-bar Implementation for Matlab - Approach 1}\label{DbarAlg_zloop}
\begin{algorithmic}[1]
\State Setup Phase. Define parameters and compute functions independent of the dynamic DN data, including:
\vspace{-2pt}
\begin{itemize}  \itemsep0pt
\item[-] Physical parameters
\item[-] Boundary parameterization and arclength function
\item[-] Current pattern matrix $J$
\item[-] Demodulate the reference data set and form the DN matrix for the reference data set $L_{\sigma_{ref}}$
\item[-] Computational grids in $k$-plane and $z$-plane
\item[-] Coefficients $\vec{c}(k),\vec{d}(k)$ for expansions of $e^{i\bar{k}\bar{z}}\vert_{\DOm}, e^{ikz}\vert_{\DOm}$ used in $\Texp$
\end{itemize}
\State Load a single frame of the measured data and demodulate
\State Compute matrix approximation to DN map, $L_{\sigma_d}$
\State Compute $\Texp$ simultaneously for all $z$ using vector operations
\State Compute the pointwise multiplication operator $T_R$ simultaneously for all $z$ using vector operations
	\ParFor{$z$ in domain} \label{zloop1}
		\State Solve D-bar equation for $\muEXPR(z,k)$
		\State $\sigma(z) \gets {(\muEXPR)}^2 (z,0)$
       \EndParFor
\end{algorithmic}
\end{algorithm}

\begin{algorithm}[h]
\caption{Fast Parallelized D-bar Implementation for Matlab - Approach 2}\label{DbarAlg_frames}
\begin{algorithmic}[1]
\State Load a batch of frames of the measured data and domain boundary points
\State Setup Phase. Define parameters and compute functions independent of the dynamic DN data, including: 
\vspace{-2pt}
\begin{itemize}  \itemsep0pt
\item[-] Physical parameters
\item[-] Boundary parameterization and arclength function
\item[-] Current pattern matrix $J$
\item[-] Demodulate the reference data set and form the DN matrix for the reference data set $L_{\sigma_{ref}}$
\item[-] Demodulate the voltage data
\item[-] Computational grids in $k$-plane and $z$-plane
\item[-] Coefficients $\vec{c}(k), \vec{d}(k)$ for expansions of $e^{i\bar{k}\bar{z}}\vert_{\DOm}, e^{ikz}\vert_{\DOm}$ used in $\Texp$
\end{itemize}
\ParFor{frames}
	\State Compute matrix approximation to DN map, $L_{\sigma_d}$
	\State Compute $\Texp$ simultaneously for all $z$ using vector operations
	\State Compute the pointwise multiplication operator $T_R$ simultaneously for all $z$ using vector operations
	\ForAll{$z$ in domain} \label{zloop2}
		\State Solve D-bar equation for $\muEXPR(z,k)$
		\State $\sigma(z) \gets {(\muEXPR)}^2 (z,0)$
	\EndFor
\EndParFor
\end{algorithmic}
\end{algorithm}

We now discuss the computational steps in detail, including numerical approximations of various functions and operations, and numerical solution of the D-bar equation. The computation of the matrix approximation to the DN map $\Lambda_{\sigma}$ can be accomplished efficiently with inner products, as explained in \cite{JenSamuBook}.  For  bipolar current patterns, as applied by ACE 1, the current pattern matrix is first transformed to an orthonormal basis $\{J^m_l\}$, $m=1,\ldots, N$ and $l=1,\ldots,L$, where $N$ is the number of linearly independent current patterns, and $L$ is the number of electrodes.  The number of linearly independent current patterns depends on the particular choice of current pattern.  For a bipolar current pattern that skips $\alpha$ electrodes between injection electrodes, $N=L-\alpha - 1$.   The voltages are transformed accordingly through the change-of-basis formula to the set denoted by $\{V^m_l\}$.  The discrete ND map $R_{\sigma_d}$ for a given data set $\sigma_d$ is approximated by
\begin{eqnarray} \label{Rmn}
R_{\sigma_d}(m,n) \approx \frac{1}{A}\sum_{l=1}^L s_l \Delta\theta_l J_l^m V_l^n,
\end{eqnarray}
where $A$ is the area of an electrode, $s_l$ is the arclength function $s(\theta)$ for the boundary evaluated at the angle $\theta_l$ corresponding to the center of the $l^{th}$ electrode,  and $\Delta\theta_l= \theta_l - \theta_{l-1}$ (where $\theta_{-1} = \theta_L)$.  Since the voltages sum to zero, we can then compute the matrix approximation $L_{\sigma_d}$ to $\Lambda_{\sigma_d}$ from $L_{\sigma_d}=(R_{\sigma_d})^{-1}$.

The scattering transform is computed for $|k|\leq R$, where $R$ is chosen empirically.  While better reconstructions can often be obtained by considering a non-uniform truncation radius $R$, we consider a disk for simplicity.  Non-uniform truncation would not result in appreciable loss of computational speed since over 95\% of the computation time is spent in the solution of the D-bar equation.  Since difference images from a reference frame are being reconstructed here, we implement the scattering transform $\Texpdif$, introduced in \cite{ChestPaper}, in which the DN map for the conductivity $\sigma=1$ is replaced by that of a reference frame $\Lambda_{\sigma_{ref}}$.  Then
\begin{equation}  \label{Texp_dif}
  \Texpdif(k) \equiv \int_{\partial\Om}e^{i\overline{k}\bar{z}}(\Lambda_{ \sigma}-\Lambda_{\sigma_{ref}})e^{ikz} ds.
\end{equation}
The functions $e^{ikz}$ and $e^{i\bar{k}\bar{z}}$ are expanded in the orthonormalized current pattern basis.
The coefficients of the expansions of the functions $e^{i\bar{k}\bar{z}}\vert_{\DOm}$ and $e^{ikz}\vert_{\DOm}$ in the orthonormalized current pattern basis vectors are computed in the setup phase of the algorithm, and are denoted by $\vec{c}(k)= [c_1(k), \ldots c_L(k)]^T$ and $\vec{d}(k)$, respectively, where
\begin{eqnarray*}
\vec{c}(k) = J^T*exp(ik\vec{z})^T \quad \mbox{and} \quad \vec{d}(k) = J^T*exp(i\bar{k}\bar{\vec{z}})^T.
\end{eqnarray*}
Then, denoting the discrete inner product over two vectors $u$ and $v$ of length $L$ by $(u,v)_L$,
\begin{eqnarray*}
\Texpdif(k) &\approx & \sum_{j=1}^{L-1} \sum^{L-1}_{m= 1}c_m(k) d_j(k) (J^j,(\Lambda_{ \sigma_d}-\Lambda_{\sigma_{ref}})J^m)_L \\
&\approx & \sum_{j=1}^{L-1} \sum^{L-1}_{m= 1}c_m(k) d_j(k) (L_{\sigma_d}(j,m)-L_{\sigma_{ref}}(j,m)) \\
\end{eqnarray*}
The fast evaluation of this formula is accomplished using inner products and vector operations.

A single grid and multigrid (2-grid) method were introduced in \cite{FIST} for the fast computation of Lippmann-Schwinger type equations that arise in D-bar methods for EIT, closely based on the method of Vainikko \cite {Vainikko}.  The convolutions are computed as FFTs, and the solution of the resulting linear system by a matrix-free method, such as GMRES.  For difference images, and in particular the data sets considered here, the frame-to-frame change in the data is sufficiently small that the method nearly always converges in one inner and one outer iteration of GMRES. In such cases, the 1-grid method is significantly faster than the 2-grid method described in \cite{FIST}.  

To apply this method, the D-bar equation \eqref{Dbar_eq} is formulated as an integral equation 
\begin{equation}\label{Fred}
  \muEXPR(z,s)=  1+\frac{1}{(2\pi)^2}\int_{|k|\leq R}\frac{\Texp(k)}{(s-k)\overline{k}}e_{-z}(k)\overline{\muEXPR(z,k)}dk_1 dk_2.
\end{equation}
To construct the computational $k$-grid, we define the square $[-s,s]^2$ where $s \geq R$, and we choose $M = 2^n$ for some positive integer $n$. The step-size for the $k$-grid is defined to be $h = 2s/(M-1)$, and the final size of the grid is then $M \times M$. We further choose  a computational $z$-grid of domain points, which need not be equally spaced.
Equation \eqref{Fred} can be written compactly as the linear system
\begin{equation} \label{LinSys}
[I - \mathcal{A}_R\; T_R( \overline{ \; \cdot \; })] \muEXPR = 1,
\end{equation}
where $T_R$ is the pointwise multiplication operator defined by
\begin{equation}
T_Rw(k) =\frac{\TexpR(k)}{4 \pi \bar{k}} e_{-z}(k) \; w(k),
\end{equation}
and the action of the operator $\mathcal{A}_R$ is given by
\begin{equation}
\mathcal{A}_R w(s)= \frac{1}{\pi} \int_{|k|<R} \frac{w(s)}{s-k} dk_1 dk_2.
\end{equation}
In our fast implementation, we compute $T_R$ simultaneously for all $z$-values in the computational grid using vector operations,  which in Matlab is more efficient than computing each $T_R$ separately inside the $z$-loop shown in step~\ref{zloop1} of Algorithm~\ref{DbarAlg_zloop} and step~\ref{zloop2} of Algorithm~\ref{DbarAlg_frames}.   The action of the operator $\mathcal{A}_R$ can be approximated by
\begin{equation}
\mathcal{A}_R w(k) \approx h^2 \mathrm{IFFT} ( \mathrm{FFT} ( \beta(k)) \cdot \mathrm{FFT}( w(k))),
\end{equation}
where $\beta(k) = (\pi k)^{-1}$ is the Green's function for the $\dbar$ operator, and $\cdot$ denotes element-wise multiplication. In Matlab, this operation can be performed efficiently using IFFTN and FFTN.  

Once the solution $\muEXPR(z,k)$ to the D-bar equation has been found, the conductivity is given by $\sigma(z) = (\muEXPR)^2 (z,0).$ Note that the D-bar equation requires the solution of the linear system \eqref{LinSys} for all $k$ values in the computational grid, even though we are ultimately only interested in $k=0$.  It is also to be noted that this method allows for the reconstruction of $ \sigma$ pointwise in $\Om$, independent from any other $z$ value, and so it is trivially parallelizable over the values in the computational $z$-grid.

For maximum computational speed, we invoked Matlab with the flag -singleCompThread, which limits Matlab to a single computational thread. This choice is compatible with the Parallel Computing Toolbox, and will restrict each parallel computation to a single core, which proved to be advantageous in the runtimes.  Several choices of numerical solvers for the linear system \eqref{LinSys} were compared for the fast implementation.  Here, we use a customized GMRES, as opposed to the built-in Matlab routine. This avoids significant overhead in calling external .m files and performing unnecessary error-checking. 

\section{Results and discussion}\label{sec:results}
To demonstrate the feasibility of clinically useful real-time reconstructions using the D-bar method, we present reconstructions from  data collected on the ACE 1 (Active Complex Electrode) EIT system at CSU.  The ACE 1 system is a bipolar current injection system with 32 active electrodes operating at a user-specified frequency up to 125 kHz.  Details of the hardware can be found in \cite{hardware}.  For the results presented here, difference images of perfusion in a cross-section of the chest of a healthy male subject sitting upright and holding his breath are presented.  360 frames of data were collected at 16 frames/s at 125 kHz and current amplitude 0.823 mA.   One frame was chosen as a reference data set and 359 difference images were computed using the fast implementation in Section \ref{sec:fast} on a uniform $k$-grid of size 16 by 16 (256 elements) and a truncation radius of $R=3.8$.  All programming was in Matlab and utilized the Parallel Computing Toolbox for the parallel solution of the D-bar equation.  

In Tables~\ref{table:runtimes_eitmac12_zloop} and \ref{table:runtimes_eitmac12_frames} we compare the performance for each version of the algorithm with various numbers of cores in parallel on a 12 core Mac Pro with two 2.66 GHz 6 core Intel Xeon processors using three different spatial grids consisting of a coarse mesh with 562 elements, a medium mesh with 1931 elements, and a fine mesh with 5,916 elements.  In Tables~\ref{table:runtimes_tx1_zloop} and \ref{table:runtimes_tx1_frames} we compare the performance on a 64 core Linux system with four 2.3 GHz 16 core  processors and 512 GB of RAM on the same three spatial grids.

Comparing all four tables, it is immediately clear that all runtimes utilizing the same number of cores were faster on the Mac Pro than on the Linux system, which is likely due in part to the difference in processor speed. It is also evident that while adding increasing numbers of parallel cores continues to improve runtimes for Algorithm~\ref{DbarAlg_frames} in Table~\ref{table:frames}, the runtimes for Algorithm~\ref{DbarAlg_zloop} in Table~\ref{table:zloop} reach maximum efficiency with a smaller number of cores, after which adding additional parallel cores actually slows the computation. This optimum number of cores increases as the $z$-mesh size increases.

It is well-known in parallel computing that the efficiency gained by adding additional parallel processors follows a ``law of diminishing returns,"  embodied by Amdahl's law \cite{Amdahl}, which states that the maximum theoretical speedup $s$ obtainable when using $n$ processors in parallel is given by
\begin{equation}\label{eq:Amdahl}
s(n) = \frac{1}{(1-p) + \frac{p}{n}},
\end{equation}
where $p$ is the proportion of the program that is parallelizable. One can see that as we increase $n$ to $\infty,$ the maximum theoretical speedup goes to $1/(1-p)$. To compute the actual speedup obtained using $n$ processors, we use $s_{actual}(n) = T(1)/ T(n)$, where $T(n)$ is the runtime with $n$ cores in parallel. 

\begin{figure}
\centering
\begin{subfigure}[h]{0.75\textwidth}
\includegraphics[width = \textwidth]{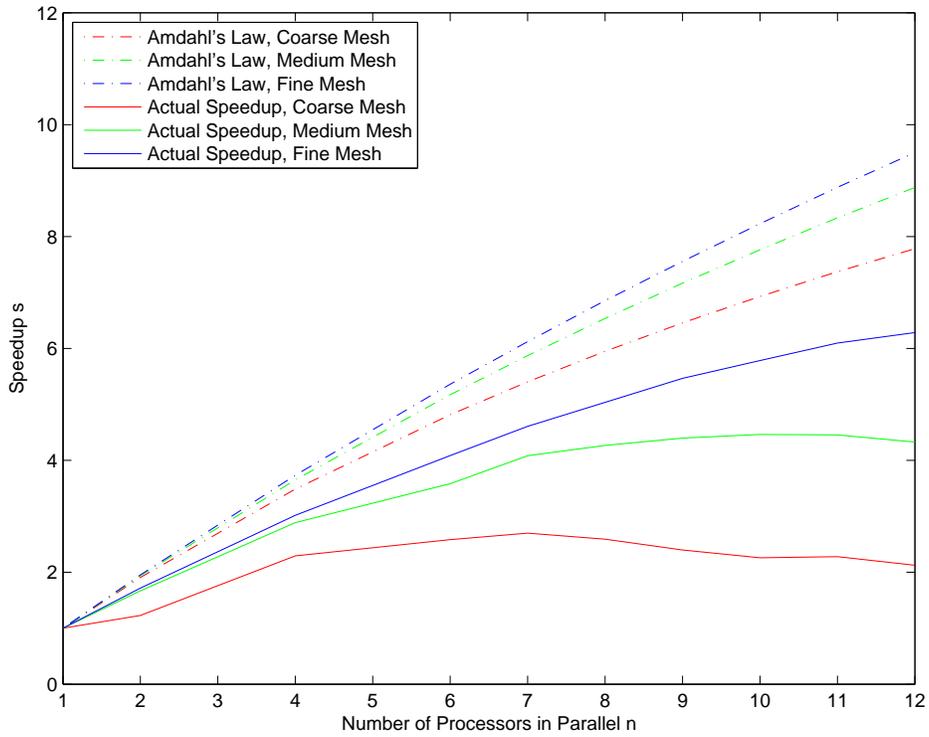}
\caption{Results with Algorithm~\ref{DbarAlg_zloop} (parallelization over mesh points)}.
\label{fig:Amdahl_zloop}
\end{subfigure}
\begin{subfigure}[h]{ 0.75\textwidth}
\includegraphics[width = \textwidth]{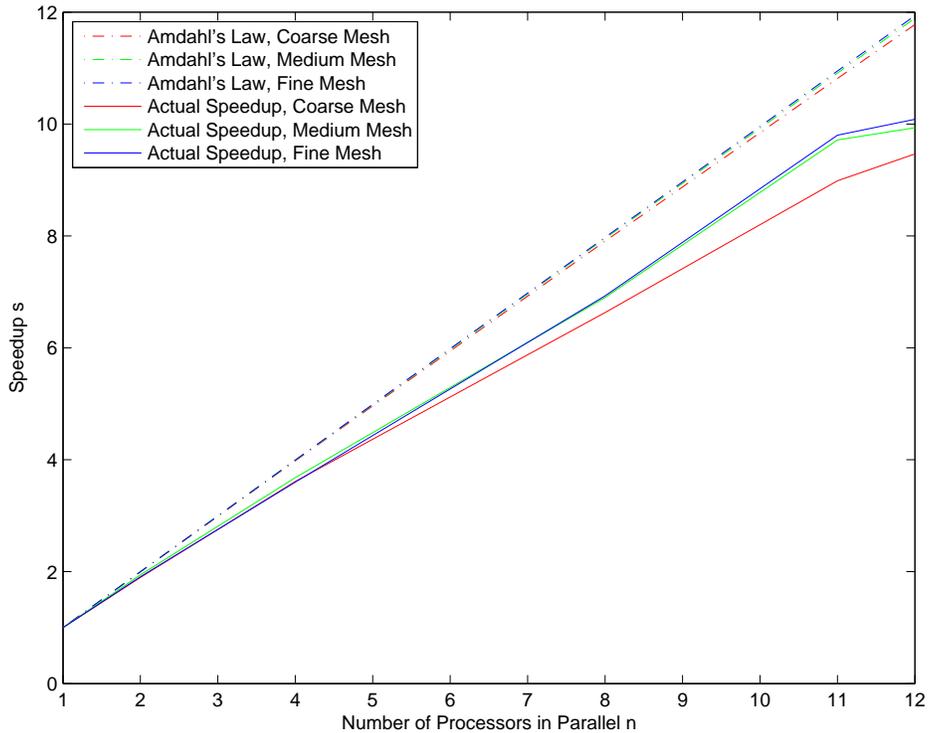}
\caption{Results with Algorithm~\ref{DbarAlg_frames} (parallelization over frames).}
\label{fig:Amdahl_frames}
\end{subfigure}
\caption{A comparison of Amdahl's law for maximum theoretical speedup (dashed lines) when using multiple cores with actual speedup obtained on  a 12 core Mac Pro with two 2.66 GHz 6 core Intel Xeon processors (solid lines) using Algorithms~\ref{DbarAlg_zloop}  and ~\ref{DbarAlg_frames}.}
\label{fig:Amdahl_eitmac12}
\end{figure}

\begin{figure}
\centering
\begin{subfigure}[h]{0.75\textwidth}
\includegraphics[width = \textwidth]{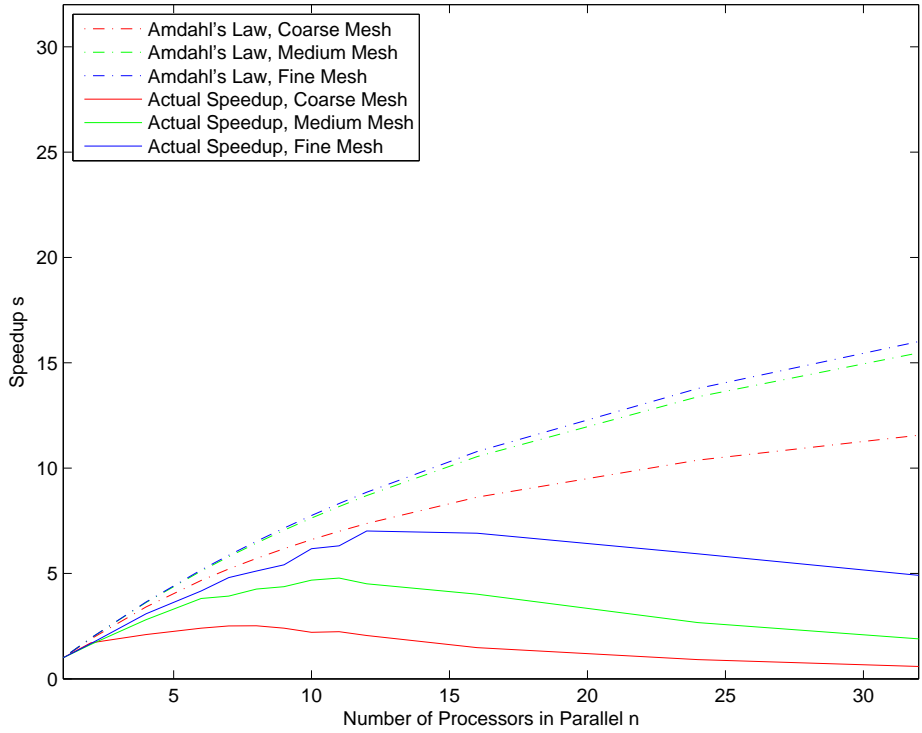}
\caption{Results with Algorithm~\ref{DbarAlg_zloop} (parallelization over mesh points).}
\label{fig:Amdahl_zloop}
\end{subfigure}
\begin{subfigure}[h]{ 0.75\textwidth}
\includegraphics[width = \textwidth]{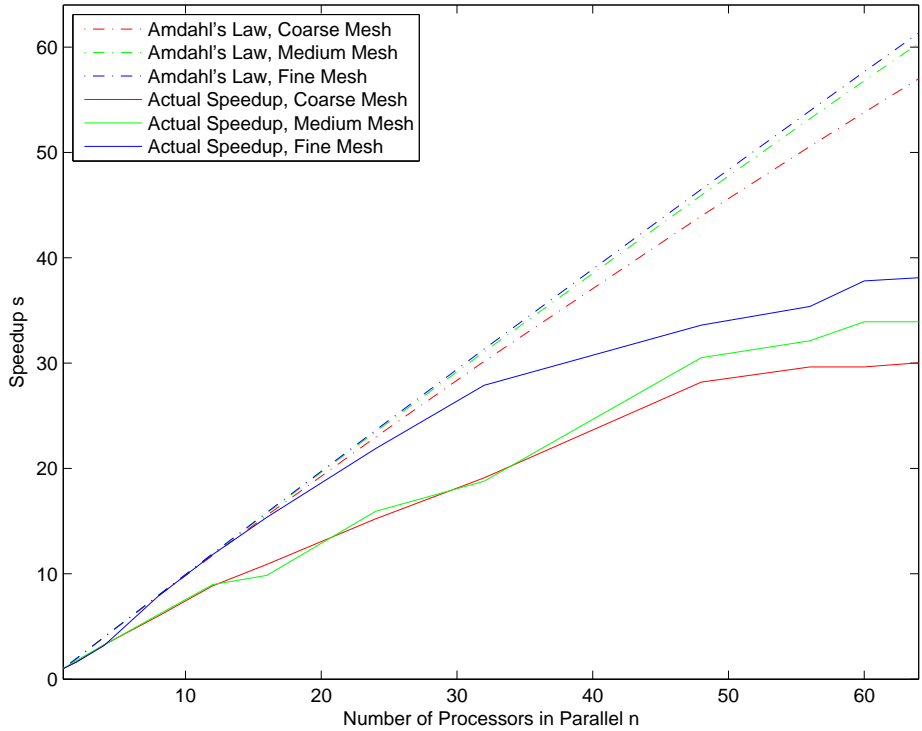}
\caption{Results with Algorithm~\ref{DbarAlg_frames} (parallelization over frames).}
\label{fig:Amdahl_frames}
\end{subfigure}
\caption{A comparison of Amdahl's law for maximum theoretical speedup (dashed lines) when using multiple cores with actual speedup obtained on   64 core Linux system with four 2.3 GHz 16 core  processors and 512 GB of RAM  (solid lines) using Algorithms~\ref{DbarAlg_zloop}  and ~\ref{DbarAlg_frames}.}
\label{fig:Amdahl_tx1}
\end{figure}

In Figure~\ref{fig:Amdahl_eitmac12}, the actual speedups obtained on the Mac Pro using both versions of the parallel algorithm are shown along with the theoretical maximum speedups predicted by \eqref{eq:Amdahl}. One can see that the results obtained using Algorithm~\ref{DbarAlg_frames} are much closer to Amdahl's ideal values than the results obtained using Algorithm~\ref{DbarAlg_zloop}. We can also see that the larger the size of the $z$-mesh, the more efficient the parallelization becomes; this is predicted by \eqref{eq:Amdahl} since increasing the number of elements in the $z$-mesh also increases $p$, the parallelizable portion of the program. In Figure~\ref{fig:Amdahl_tx1} the same plots are included for the 64 core Linux system.  It is evident that the additional cores do not provide speedup for Algorithm~\ref{DbarAlg_zloop}, and the divergence from the theoretically predicted speedup by Amdahl's Law increases with the number of cores, but Amdahl's Law also predicts the leveling off of speedup at around 20 to 25 cores; in fact, this happens at around 10 to 12 cores for parallelization over mesh values (Algorithm~\ref{DbarAlg_zloop}).  However, for Algorithm~\ref{DbarAlg_frames}, the speedup continues for all 64 cores, although it is clearly starting to level out at around 60 cores.

\begin{subtables}
\begin{table}[p]
\centering
\rotatebox{90}{
\begin{minipage}{1.17\linewidth}
\centering
\begin{tabular*}{1.17\linewidth}{|c|r|r|r|r|r|r|r|r|r|}
\toprule    
      & \multicolumn{3}{c|}{Coarse grid (562 $z$-values)}   &
       \multicolumn{3}{|c|}{Medium grid (1931 $z$-values)}  & \multicolumn{3}{c|}{Fine grid (5,916 $z$-values)}  \\  \midrule
Cores &  Total runtime (s) & Loop runtime (s) &s/frame  & Total runtime (s) & Loop runtime (s) & s/frame& Total runtime (s) & Loop runtime (s) & s/frame\\
        \midrule
1	&60.1527	&57.1906	&0.1676	&200.1049	&193.7054	&0.5574	&2449.0000	&2397.3000	&6.8217	\\
2	&49.0169	&45.0331	&0.1365	&119.8318	&113.8120	&0.3338	&1450.8000	&1384.0000  &4.0412 	\\
4	&26.2268	&23.1675	&0.0731	&69.2834	&62.6818	&0.1930	&793.1958	&719.3565	&2.2095	\\
6	&23.3147	&20.3212	&0.0649	&55.8466	&47.2439	&0.1556	&581.6291	&508.4273	&1.6201	\\
7	&22.2798	&19.2730	&\hl{0.0621}	&49.0180	&41.6290	&0.1365	&516.3475	&442.8791	&1.4383	\\
8	&23.2349	&20.2471	&0.0647	&46.8985	&38.7851	&0.1306	&465.2694	&393.8668	&1.2960	\\
9	&25.0864	&22.0817	&0.0699	&45.4834	&37.2894	&0.1267	&423.6890	&352.4819	&1.1802	\\
10	&26.6466	&23.6232	&0.0742	&44.8317	&36.8733	&\hl{0.1249}	&393.7915	&322.4939	&1.0969	\\
11	&26.3834	&23.3987	&0.0735	&44.8993	&36.9119	&0.1251	&371.5604	&300.3748	&1.0350	\\
12	&28.2897	&25.3129	&0.0788	&46.2567	&38.1605	&0.1288	&355.7963	&283.7342	&\hl{0.9911}	\\
       \bottomrule
\end{tabular*}
\caption{Runtimes in seconds for Algorithm~\ref{DbarAlg_zloop} (parallelization over mesh points) over 359 frames on a 12 core Mac Pro with two 2.66 GHz 6 core Intel Xeon processors. Here, the loop runtimes refer to the runtime for the parallelized loop over $z$-values. The fastest per-frame runtime for each mesh has been highlighted.}
\label{table:runtimes_eitmac12_zloop}
\end{minipage} } 
\qquad
\rotatebox{90}{
\begin{minipage}{1.17\linewidth}
\centering
\begin{tabular*}{1.17\linewidth}{|c|r|r|r|r|r|r|r|r|r|}
\toprule         
      & \multicolumn{3}{c|}{Coarse grid (562 $z$-values)}   &
       \multicolumn{3}{|c|}{Medium grid (1931 $z$-values)}  & \multicolumn{3}{c|}{Fine grid (5,916 $z$-values)}  \\  \midrule
Cores &  Total runtime (s) & Loop runtime (s) &s/frame  & Total runtime (s) & Loop runtime (s) & s/frame& Total runtime (s) & Loop runtime (s) & s/frame\\
       \midrule
1	&70.6320	&66.6039	&0.1967	&239.3496	&231.0972	&0.6667	&767.5243	&742.7878	&2.1380	\\
2	&41.4227	&38.4427	&0.1154	&146.6424	&137.5511	&0.4085	&460.4769	&435.7931	&1.2827	\\
4	&33.5734	&28.9639	&0.0935	&85.1865	&75.0827	&0.2373	&248.8144	&229.7575	&0.6931	\\
6	&29.3667	&24.9876	&0.0818	&62.8336	&55.3416	&0.1750	&184.0372	&164.8297	&0.5126	\\
7	&28.1477	&24.1555	&0.0784	&60.9286	&51.3155	&0.1697	&159.8174	&142.7302	&0.4452	\\
8	&28.0048	&24.1189	&\hl{0.0780}	&56.1551	&46.3121	&0.1564	&150.0559	&127.9969	&0.4180	\\
11	&31.5985	&26.9765	&0.0880	&50.0699	&42.1989	&\hl{0.1395}	&121.6363	&99.6636	&0.3388	\\
12	&34.3986	&30.1875	&0.0958	&53.0852	&40.8174	&0.1479	&109.3473	&91.8146	&\hl{0.3046}	\\
16	&47.7978	&43.7131	&0.1331	&59.5714	&50.6782	&0.1659	&111.0096	&83.2492	&0.3092	\\
24	&77.3859	&72.0929	&0.2156	&89.5173	&80.3598	&0.2494	&129.2778	&104.2357	&0.3601	\\
32	&120.9989	&116.1762	&0.3370	&125.9324	&115.8191	&0.3508	&156.2232	&131.1057	&0.4352	\\
       \bottomrule
\end{tabular*} 
\caption{Runtimes in seconds for Algorithm~\ref{DbarAlg_zloop} (parallelization over mesh points) over 359 frames on a 64 core Linux system with four 2.3 GHz 16 core  processors and 512 GB of RAM. Here, the loop runtimes refer to the runtime for the parallelized loop over the frames. The fastest per-frame runtime for each mesh has been highlighted.}
\label{table:runtimes_tx1_zloop}
\end{minipage} }
\end{table}
\label{table:zloop}
\end{subtables}

\begin{subtables}
\begin{table}[p]
\centering
\rotatebox{90}{
\begin{minipage}{1.17\linewidth}
\centering
\begin{tabular*}{1.17\linewidth}{|c|r|r|r|r|r|r|r|r|r|}
\toprule    
      & \multicolumn{3}{c|}{Coarse grid (562 $z$-values)}   &
       \multicolumn{3}{|c|}{Medium grid (1931 $z$-values)}  & \multicolumn{3}{c|}{Fine grid (5,916 $z$-values)}  \\  \midrule
Cores &  Total runtime (s) & Loop runtime (s) &s/frame  & Total runtime (s) & Loop runtime (s) & s/frame& Total runtime (s) & Loop runtime (s) & s/frame\\
        \midrule
1	&63.5468	&63.4370	&0.1770	&206.1180	&205.9474	&0.5741	&623.0107	&622.7078	&1.7354 	\\
2	&33.5405	&33.4232	&0.0934	&106.1961	&106.0467	&0.2958	&326.5889	&326.2849	&0.9097	\\
4	&17.5768	&17.4568	&0.0490	&55.9864	&55.8287	&0.1560	&173.0110	&172.6940	&0.4819	\\
8	&9.5713	&9.4459	&0.0267	&29.8543	&29.7101	&0.0832	&89.9303	&89.6408	&0.2505	\\
11	&7.0793	&6.9545	&0.0197	&21.2018	&21.0424	&0.0591	&63.5926	&63.2956	&0.1771	\\	
12	&6.7303	&6.6036	&\hl{0.0187}	&20.7440	&20.5766	&\hl{0.0578}	&61.7965	&61.4933	&\hl{0.1721}	\\
        \bottomrule
\end{tabular*}
 \caption{Runtimes in seconds for Algorithm~\ref{DbarAlg_frames} (parallelization over frames) over 359 frames on a 12 core Mac Pro with two 2.66 GHz 6 core Intel Xeon processors. Here, the loop runtimes refer to the runtime for the parallelized loop over the frames. The fastest per-frame runtime for each mesh has been highlighted.}
 \label{table:runtimes_eitmac12_frames}
\end{minipage} } 
\qquad
\rotatebox{90}{
\begin{minipage}{1.17\linewidth}
\centering
\begin{tabular*}{1.17\linewidth}{|c|r|r|r|r|r|r|r|r|r|}
\toprule         
      & \multicolumn{3}{c|}{Coarse grid (562 $z$-values)}   &
       \multicolumn{3}{|c|}{Medium grid (1931 $z$-values)}  & \multicolumn{3}{c|}{Fine grid (5,916 $z$-values)}  \\  \midrule
Cores &  Total runtime (s) & Loop runtime (s) &s/frame  & Total runtime (s) & Loop runtime (s) & s/frame& Total runtime (s) & Loop runtime (s) & s/frame\\
       \midrule
1	&83.0041	&82.8420	&0.2312	&261.7762	&261.5257	&0.7292	&837.1351	&836.5582	&2.3319	\\
2	&49.8535	&49.5505	&0.1389	&145.1565	&145.9066	&0.4043	&510.6861	&510.2197	&1.4225	\\
4	&25.3656	&25.1796	&0.0707	&80.8684	&80.6215	&0.2253	&262.7161	&262.2612	&0.7318	\\
8	&13.8435	&13.6768	&0.0386	&42.7368	&42.4592	&0.1190	&106.2459	&105.8588	&0.2959	\\
12	&9.3821	&9.1925	&0.0261	&29.1641	&28.8975	&0.0812	&70.7996	&70.4057	&0.1972	\\
16	&7.6005	&7.3360	&0.0212	&26.5531	&26.3019	&0.0740	&54.4559	&54.0610	&0.1517	\\	
24	&5.4611	&5.2730	&0.0152	&16.4534	&16.1881	&0.0458	&38.2854	&37.8909	&0.1066	\\
32	&4.3561	&4.1731	&0.0121	&13.9393	&13.6828	&0.0388	&30.0064	&29.5878	&0.0836	\\
48	&2.9294	&2.7678	&0.0082	&8.5636	&8.3324	&0.0239	&24.9113	&24.4925	&0.0694	\\
60	&2.7924	&2.6435	&0.0078	&7.7249	&7.4911	&0.0215	&22.1673	&21.7229	&0.0617	\\
64	&2.7781	&2.6219	&\hl{0.0077} &7.7157	&7.4769	&\hl{0.0215} &21.9644	&21.5439	&\hl{0.0612}	\\
        \bottomrule
\end{tabular*}
 \caption{Runtimes in seconds for Algorithm~\ref{DbarAlg_frames} (parallelization over frames) over 359 frames on a 64 core Linux system with four 2.3 GHz 16 core  processors and 512 GB of RAM. Here, the loop runtimes refer to the runtime for the parallelized loop over the frames. The fastest per-frame runtime for each mesh has been highlighted.}
 \label{table:runtimes_tx1_frames}
\end{minipage} }
\end{table} \label{table:frames}
\end{subtables}

\begin{figure}[h] 
\begin{center}
\includegraphics[width=\textwidth]{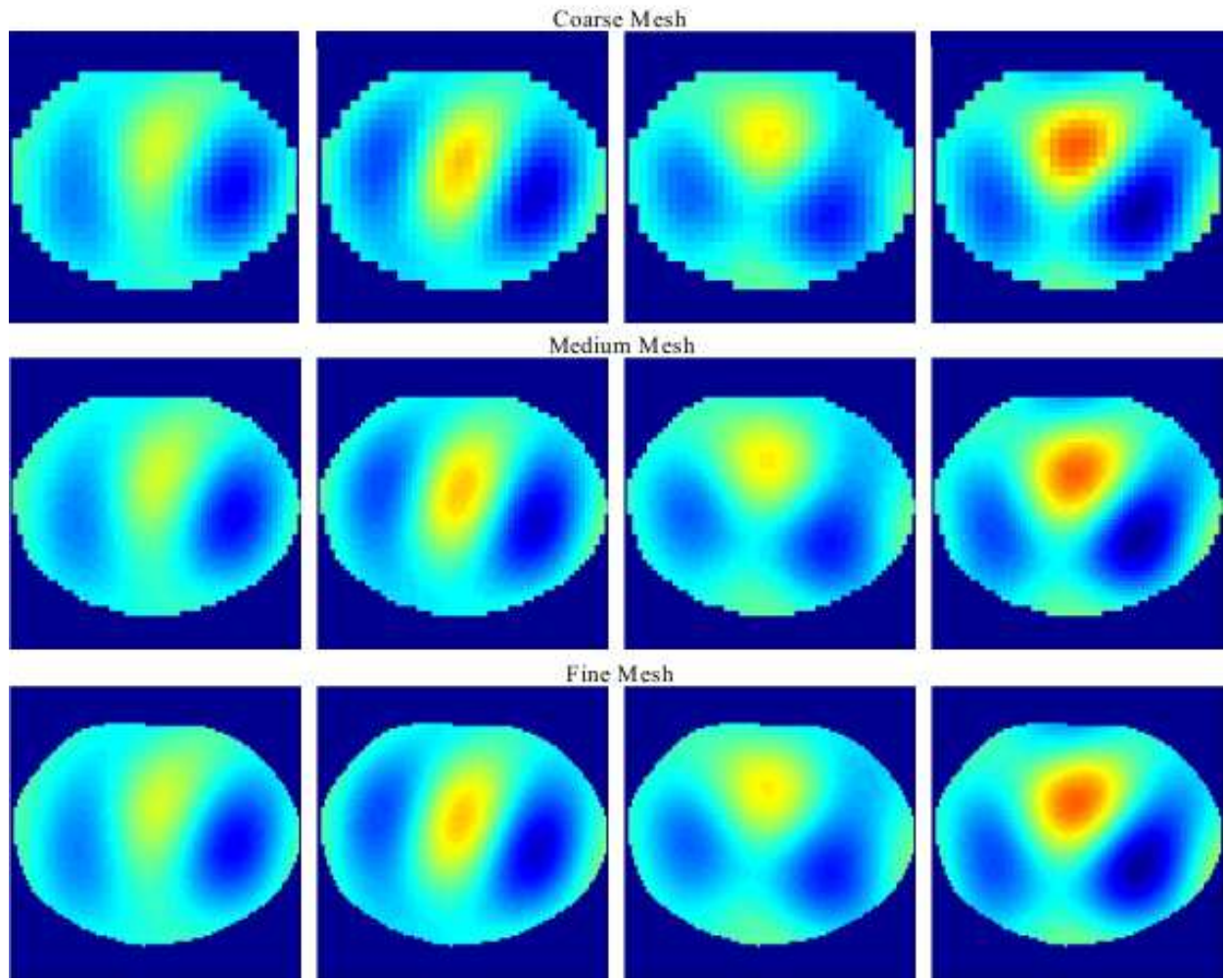} 
\end{center}
\caption{Reconstructions on the three $z$-meshes  timed in this work of four frames in the sequence of 360 frames showing changes due to perfusion in the chest of a healthy human subject.  The heart is at the top, and red represents high conductivity and blue low conductivity with respect to the reference frame.  The images are displayed on the same scale.  }
\label{recon}
\end{figure}

Figure~\ref{recon} contains four frames in the reconstruction of the human chest data displayed in the three $z$-meshes to illustrate the resolution provided by these mesh choices.    The figure depicts changes due to perfusion.  The heart is at the top, and red represents high conductivity and blue low conductivity.  The images are displayed on the same scale.   While the topmost mesh is quite coarse, the lungs and heart are clearly visible, and changes are evident.  The medium mesh is significantly better, and probably provides the best compromise for real-time imaging, but this implementation comes with the price of an approximately 1 second delay.  The very fine mesh is included to illustrate highly desirable spatial resolution and the associated runtimes, given in Tables~\ref{table:zloop} and \ref{table:frames}.

From the clinical perspective, difference images are sufficient for some applications, such as real-time detection of a pneumothorax or atelectasis.  For other applications, such as distinguishing between blood, water, and mucus in the lung, absolute images may be required, which may require longer computation times.  Some applications may also require finer spatial resolution than that presented here.  However, the data is preserved for off-line or delayed reconstruction, and improvements such as the use of additional cores in parallel, faster FFTs, and faster processors are likely to yield improved runtimes in the future.

\section{Conclusions}
The results presented here show for the first time the D-bar method applied to human chest data collected on a pairwise current injection system.  Conductivity changes due to perfusion are clearly visible in the images.   The fast implementation demonstrates the clinical potential of the D-bar algorithm as a reconstruction algorithm for real-time bedside imaging.

\section{Acknowledgments} The project described was supported by Award Number 1R21EB016869-01A1 from the National Institute Of Biomedical Imaging And Bioengineering.  The content is solely the responsibility of the authors and does not necessarily represent the official view of the National Institute Of Biomedical Imaging And Bioengineering or the National Institutes of Health.

\end{document}